\documentclass[11pt]{article}

\usepackage{amsmath}
\usepackage{amsfonts}
\usepackage{amssymb}
\usepackage{amsthm}
\usepackage{graphics}
\usepackage{amscd}
\usepackage{graphicx,epsfig}
\usepackage{color}
\usepackage[all]{xy}
\usepackage{url}

\DeclareMathOperator{\Div}{div}

\renewcommand{\epsilon}{\varepsilon}

\newcommand{\boD}{\mathcal{D}}

\newcommand{\boC}{\mathcal{C}}

\newcommand{\boS}{\mathcal{S}}

\newcommand{\boN}{\mathcal{N}}

\newcommand{\R}{\mathbb{R}}

\renewcommand{\S}{\mathbb{S}}
\newcommand{\N}{\mathbb{N}}
\newcommand{\der}[2]{\dfrac{\partial #1}{\partial #2}}
\newcommand{\dd}{\mathrm{d}}

\newcommand{\Ome}{\Omega}

\newtheorem{thm}{Theorem}

\newtheorem{prop}[thm]{Proposition}
\newtheorem{lem}[thm]{Lemma}

\renewcommand{\phi}{\varphi}
\newcommand{\dis}{\displaystyle}

\newtheorem*{thm*}{Theorem}

\theoremstyle{remark}

\newtheorem*{rem*}{Remark}

\newcounter{remark}

\newcounter{case}

\newcounter{construction}

\newcounter{fact}

\title{Some uniqueness results for constant mean curvature graphs}
\author{Laurent Mazet}
\date{}

\begin{document}

\maketitle

\begin{abstract}
The aim of this paper is to give two uniqueness results for the Dirichlet
problem associated to the constant mean curvature equation. We study constant
mean curvature graphs over strips of $\R^2$. The proofs are based on height
estimates and the study of the asymptotic behaviour of solutions to the
Dirichlet problem.
\end{abstract}



\section*{Introduction}
The surfaces with constant mean curvature are the mathematical modelling of
soap films. These surfaces appear as the interfaces in isoperimetric
problems. There exist different points of view on constant mean curvature
surfaces, one is to consider them as graphs.

Let $\Ome$ be a domain of $\R^2$. The graph of a function $u$ over $\Ome$ has
constant mean curvature $H>0$ if it satisfies the following partial
differential equation:
\begin{equation*}
\label{cmc}
\Div\left(\frac{\nabla u}{\sqrt{1+|\nabla u|^2}}\right)=2H
\tag{CMC}
\end{equation*}
The graph of such a solution is called a $H$-graph and has a upward pointing
mean curvature vector.

Thanks to the work of J.~Serrin \cite{Se1,Se2} and J.~Spruck \cite{Sp}, we can
build a lot of $H$-graphs over bounded domains of $\R^2$. Over unbounded
domains, the Dirichlet problem associated to \eqref{cmc} is more
complicated. Existence results are known and due to P.~Collin \cite{Co} and
R.~L\'opez \cite{Lo1,Lo2}. The question we ask in this paper is the uniqueness
of solutions for the boundary data that these authors study.

In \cite{Lo1,Lo2}, R.~L\'opez proves existence results for vanishing boundary
data. In the case of bounded boundary data, the uniqueness is known. In
\cite{Co} and \cite{Lo1}, the authors build solutions of \eqref{cmc} over
strips with boundary data that can be unbounded. In this paper, we prove that,
for these boundary data, we have uniqueness (see Theorems \ref{uniciteconvex}
and \ref{unicitecircle}).


There are two major steps to prove these results. First, if there are two
solutions for a same boundary data, the difference between these solutions can
not stay bounded. Thus, this gives us an information on the asymptotic
behaviour of the boundary data. The second step consists in seeing the
consequences of this asymptotic behaviour on the
asymptotic behaviour of a solution. In this second step, we use the notion of
line of divergence that the author has defined in \cite{Ma1}.

The paper is devided as follows. In the first section, we give the two
existence results of P.~Collin and R.~L\'opez that we are interested in. We
also precise some definitions and notations. In the second part, we give a set
of upper and under-bounds for the $H$-graph that we study. These results are
important to prepare the proofs of our uniqueness results. Section
\ref{proofcollin} is devoted to the proof of the uniqueness of Collin's
solutions. In Section \ref{prooflopez}, we prove the uniqueness of L\'opez's
solutions. The two proofs are very similar.

\section{Two existence results}
In this first section, we recall two existence results for Dirichlet problem;
these results were proved by P.~Collin in \cite{Co} and R.~L\'opez in
\cite{Lo1}. 

In both cases, the author studies the Dirichlet problem associated to the
constant mean curvature equation \eqref{cmc} on a strip
$\Ome= \R\times (-l,l)$ of width $2l$. It was proven in \cite{Lo2}, that the
width $2l$ needs to be less than $1/H$ for having a solution.

\subsection{The existence results}

Let $f:\R\rightarrow \R$ be a continuous real function. We define $\phi_f$ on
$\partial \Ome$  by $\phi_f(x,\pm l)=f(x)$. P.~Collin and R.~L\'opez have
looked for a solution $u$ of the constant mean curvature equation
\begin{equation*}
\Div\left(\frac{\nabla u}{\sqrt{1+|\nabla u|^2}}\right)=2H
\tag{CMC}
\end{equation*}
such that $u|_{\partial\Ome}=\phi_f$.

The result of P.~Collin concerns the limiting case $2l=1/H$.
\begin{thm*}[P.~Collin, \cite{Co}]\label{collin}
Let $f:\R\rightarrow\R$ be a convex continuous function. There exists a
solution $u$ of \eqref{cmc} on $\Ome=\R\times (-1/(2H),1/(2H))$ which takes
$\phi_f$ as boundary value. 
\end{thm*}
This result was independently proved by A.~N.~Wang for the convex function
$x\mapsto x^2$ (see \cite{Wa}).

Let $\Ome$ be a domain in $\R^2$, we say that $\Ome$ satisfies a
\emph{uniform exterior $R$-circle condition} if at each point
$p\in\partial\Ome$ there exists a disk $D$ with radius $R$ such that
$\overline{D}\cap \overline{\Ome}=\{p\}$. This tells us that a circle of radius
$R$ can "roll" outside $\Ome$ along $\partial\Ome$ touching each point of
$\partial\Ome$ along its deplacement.

Let $f:\R\rightarrow\R$ be a countinuous function. $f$ satisfies a
\emph{uniform $R$-circle under condition} if the domain $\{(x,y)\in\R^2\,|\,
  y\ge f(x)\}$ satisfies a uniform exterior $R$-circle condition. It says that
a circle of radius $R$ can "roll" under the graph of $f$ touching each point
of the graph along its deplacement.

The result of R.~L\'opez deals with the case where $2l<1/H$. We give it now
using some notations that we shall introduce in the following subsection.
\begin{thm*}[R.~L\'opez, \cite{Lo1}]\label{lopez}
Let $f:\R\rightarrow\R$ be a continuous function and $t\in\R_+^*$. We assume
that $f$ satisfies a uniform $\rho_t(H)$-circle under condition. Then, there
exists a solution $u$ of \eqref{cmc} on the strip
$\Ome=\R\times(-h_t(H),h_t(H))$ which takes $\phi_f$ as boundary value. 
\end{thm*}

The technique used by P.~Collin and R.~Lopez is the Perron technique. They
build their solutions as the supremum of under-solutions. The difficulty is to
have good barrier functions to ensure the boundary value.

The aim of this paper is to prove that the solutions build by P.~Collin and
R.~L\'opez are unique for the boundary data $\phi_f$.

\begin{thm*}[Theorems \ref{uniciteconvex} and \ref{unicitecircle}]
Let $f:\R\rightarrow\R$ be a continuous function. We consider one of the two
cases below:
\begin{enumerate}
\item $\Ome=\R\times (-\frac{1}{2H},\frac{1}{2H})$ and $f$ is convex.
\item $\Ome=\R \times (-h_t(H),h_t(H))$ and $f$ satisfies a uniform
  $\rho_t(H)$-circle under condition. 
\end{enumerate}
Let $u$ and $v$ be two solutions of \eqref{cmc} on $\Ome$ with $\phi_f$ as
boundary value. Then $u=v$.
\end{thm*}

\subsection{The $1$-parameter family of nodoids}

In this subsection, we recall many classical results on nodoids, the
reader can refer to \cite{De, Ee, Lo1} for more explanations.

The constant mean curvature surfaces of revolution are the Delaunay
surfaces. This set of surfaces splits into two $1$-parameter families. One is
composed of embedded surfaces which are called unduloids. When the parameter
moves, the family of the unduloids deforms the cylinder of radius $1/(2H)$
into a stack of tangent sphere of radius $1/H$.

The second $1$-parameter family is composed of non-embedded surfaces: these
are the nodoids. The interest in considering the nodoids is that, since they
have self-intersections, each one contains a piece that looks like a
catenoidal neck with mean curvature vector pointing outward.

Let us recall the construction of nodoids and fix notations. Let $r(u)$
be a positive smooth function defined on an open interval $I$ and consider the
surface of revolution parametrized by:
$$
X(u,\theta)=\left(r(u)\cos(\theta),r(u)\sin(\theta),u\right)
$$
We fix the normal to the surface to be 
$$
N(u,\theta)=\frac{1}{\sqrt{1+r'^2}}\left(\cos(\theta),\sin(\theta),-r'\right)
$$

The surface has constant mean curvature $H$ if $r$ satisfies:
$$
2H=-\frac{1}{r\sqrt{1+r'^2}}+\frac{r''}{(1+r'^2)^{3/2}}
$$
After a first integration, this equation implies that there exists $c\in\R$
such that:
\begin{equation}\label{eq1}
Hr^2=-\frac{r}{\sqrt{1+r'^2}}+c
\end{equation}

Since $Hr^2$ is positive, $c$ needs to be positive. When $c>0$, there exist
$h$, $\rho$ and $r:[-h,h]\rightarrow [0,\rho]$ a solution to \eqref{eq1}
such that $r$ is even and the initial value $r(0)=t>0$ is the minimum 
of $r$. Besides $r(h)=\rho$ and $r'(h)=+\infty$, then $H\rho^2=c$. The
associated surface $X$ is a nodoid.

For $u=0$, we have $Ht^2+t=c$ then:
$$
t=\frac{-1+\sqrt{1+4Hc}}{2H}
$$
$t$ is an increasing function of $c$ with $t=0$ for $c=0$ and
$\lim_{c\rightarrow +\infty} t=+\infty$. In the following, $t$ is then used as
parameter for the family of nodoids. We have:
$$
\rho=\rho_t(H)=\sqrt{\frac{Ht^2+t}{H}}
$$

Moreover, we have:
$$
h=h_t(H)=\int_t^\rho\frac{H(\rho^2-x^2)}{\sqrt{x^2-H^2(\rho^2-x^2)^2}}\dd x
$$
\begin{figure}
\begin{center}
\resizebox{0.6\linewidth}{!}{\input{figuniciteH5.pstex_t}}
\caption{\label{nodoid2}}
\end{center}
\end{figure}
We can summarize the properties in the following propostion and in Figure
\ref{nodoid2}. 

\begin{prop}
There exists a $1$-parameter family of nodoids $\{\boN_t,t>0\}$ with constant
mean curvature $H$ given by the rotation of a curve $\gamma_t$ around the
$z$-axis and with the following properties:
\begin{enumerate}
\item The curve $\gamma_t$ is a graph on $[h_t(H),h_t(H)]$ of an even
  function.
\item The curve $\gamma_t$ has horizontal tangents  at $\pm h_t(H)$. Then
  $\boN_t$ is included in the slab $\boS_t$ : $|z|\le h_t(H)$ and is tangent
  to it.
\item The mean curvature vector points outside the bounded domain determined
  by $\boN_t$ in the slab $\boS_t$.
\item The circle $C_t$ of $\boN_t$ with smallest radius is given by
  $x^2+y^2=t^2,\ z=0$.
\item The function $h_t(H)$ is strictly increasing  on $t$ and 
\begin{align*}
\lim_{t\rightarrow 0}h_t(H)&=0& \lim_{t\rightarrow
  +\infty}h_t(H)&=\frac{1}{2H} 
\end{align*} 
\item The function $\rho_t(H)$ is strictly increasing and 
\begin{align*}
\lim_{t\rightarrow 0}\rho_t(H)&=0& \lim_{t\rightarrow +\infty}
\rho_t(H)&=+\infty& \lim_{t\rightarrow +\infty}\rho_t(H)-t&=\frac{1}{2H}
\end{align*}
\end{enumerate}
\end{prop}

The two limits of $h_t(H)$ and $\rho_t(H)$ as $t\rightarrow +\infty$ allow us
to consider P.~Collin result as a limiting case of R.~L\'opez
theorem. Actually, when $R$ goes to $+\infty$, the uniform $R$-circle under
condition for $f$ becomes the convexity since the circle becomes a
straight-line. 

\section{Preliminaries}

\subsection{The maximal and minimal solutions}

To prove the uniqueness of the solutions built by P.~Collin and R.~L\'opez, we
need some control on solutions of these Dirichlet problem. We have a first
result.

\begin{lem}\label{supminimal}
Let $f:\R\rightarrow \R$ be a continuous function. On $\Ome=\R\times (-l,l)$,
there exists a solution $w$ of the minimal surface equation
\begin{equation*}\label{esm}
\Div\left(\frac{\nabla w}{\sqrt{1+|\nabla w|^2}}\right)=0
\tag{MSE}
\end{equation*}
with $w=\phi_f$ on the boundary of $\Ome$. Besides, we have $w\ge u$ for every
solution $u$ of \eqref{cmc} on $\Ome$ with $u=\phi_f$ on $\partial\Ome$.
\end{lem}

\begin{proof}
Let us consider $n\in \N^*$. Because of a result by H.~Jenkins and J.~Serrin
\cite{JS}, if $n$ is big enough, there exists two solutions $w_n^+$ and
$w_n^-$ of \eqref{esm} on $(-n,n)\times (-l,l)$ with $w_n^\pm=\phi_f$ on
$(-n,n)\times\{-l,l\}$ and $w_n^\pm=\pm\infty$ on $\{-n,n\}\times (-l,l)$. 

By maximum principle, for every $n$ and $m$, we have $w_n^+\ge w_m^-$ and
$(w_n^+)$ is a decreasing sequence. This implies that $(w_n^+)$ converges to a
solution $w$ of \eqref{esm} on $\Ome$ with $\phi_f$ as boundary value.

Let us consider now a solution $u$ of \eqref{cmc} on $\Ome$ with $\phi_f$ as
boundary value. By maximum principle, for every $n$, we have $w_n^+\ge u$ ;
then in the limit, $w\ge u$.
\end{proof}

This lemma gives an upper-bound to a solution of \eqref{cmc} without any
hypothesis on the function $f$. To get an under-bound, we need such
hypotheses.

Let us consider the function $c$ which is defined on $\Ome
=\R\times(-1/(2H),1/(2H))$ by:
$$
c(x,y)=-\frac{1}{\cos\theta}\sqrt{\frac{1}{4H^2}-y^2}+(x-x_0)\tan\theta+ z_0 
$$
$c$ is a solution of \eqref{cmc}: its graph is in fact the half-cylinder with
the two straight-lines of equation $z=(x-x_0)\tan\theta +z_0$ over
$\partial\Ome$ as boundary.
\begin{lem}\label{infconvex}
Let $f:\R\rightarrow\R$ be a convex function and $u$ a solution of \eqref{cmc}
on $\Ome=\R\times(-1/(2H),1/(2H))$ with $\phi_f$ as boundary value. Let $x_0$
be in $\R$ and $z=(x-x_0)\tan\theta_0 +f(x_0)$ be a straight-line
which is below 
the graph of $f$ (such a line exists because of the convexity). Let $c$ denote
the half-cylinder associated to this line. Then we have $u\ge c$ on
$\Ome$. 
\end{lem}

\begin{proof}
Let $h$ denote the function defined on $\Ome$ by $h(x,y)=(x-x_0)\tan\theta_0
+f(x_0)$. On the boundary $u\ge h$.

If the function $f$ is affine, \textit{i.e.} $f(x)=(x-x_0)\tan\theta_0
+f(x_0)$, it is known that $c$ is the only constant mean curvature extension
for $\phi_f$ (see Theorem 8 in \cite{Ma2}), then $u=c$.

If the function $f$ is non-affine, the set of $\theta$ such that there exists
$x_1\in\R$ with $z=(x-x_1)\tan\theta +f(x_1)$ is below the graph of $f$ is an
interval $I\subset \R$. We assume that $\theta_0$ is in the interior of
this interval. If $\theta_0$ is an end-point of this interval, the property is
proved by continuity.

Since $\theta_0$ is in the interior of $I$, there exist $x_1< x_0< x_2$ and
$\theta_1<\theta_0<\theta_2$ such that $(x-x_1)\tan\theta_1 +f(x_1)\le f$ and
$(x-x_2)\tan\theta_2 +f(x_2)\le f$. By Proposition 3 in \cite{Ma2}, there
exists $K\in\R_+$ such that:
\begin{gather*}
u(x,y)\ge (x-x_1)\tan\theta_1 +f(x_1)-K\\
u(x,y)\ge (x-x_2)\tan\theta_2 +f(x_2)-K
\end{gather*}
Since $\theta_1<\theta_0<\theta_2$, these two equations imply that $u(x,y)\ge
h(x,y)$ if $|x|$ is big enough. We have $h\ge c$ on $\Ome$ (we recall that
$c$ is the half-cylinder associated to $z=(x-x_0)\tan\theta_0 +f(x_0)$); then
$u\ge c$ on $\partial\Ome$ and outside a compact of
$\overline{\Ome}$. Finally, by maximum principle, $u\ge c$ in $\Ome$.
\end{proof}

In the case of L\'opez solutions, we get the following under-bound.
\begin{lem}\label{infcircle}
Let $f:\R\rightarrow\R$ a continuous function that satisfies a uniform
$\rho_t(H)$-circle under condition. Let $x$ be in $\R$ and $\boC$ a circle of
radius $\rho_t(H)$ that established the uniform $\rho_t(H)$-circle under
condition at the point $(x,f(x))$. Let $u$ be a solution of \eqref{cmc} on
$\Ome=\R\times (-h_t(H),h_t(H))$ with $\phi_f$ as boundary value. Then, the
graph of $u$ is above the nodoid $\boN_t$ which have horizontal axis and is
bounded by the two parallel circles $\boC$ in the vertical plane $y=-h_t(H)$
and $y=h_t(H)$.
\end{lem}
\begin{proof}
Let $e_z$ denote the vertical unit vector $(0,0,1)$. For $s$ in
$\R$, let us translate by $se_z$ the nodoid $\boN_t$ bounded by
the two parallel circles $\boC$. For enough negative $s$, $\boN_t+se_z$ is
below the graph of $u$. Let $s$ grow until the first contact. The mean
curvature of the graph is upward pointing and the mean curvature of $\boN_t$
points outside. So by maximum principle, the first contact can not be an
interior point. Then, because of the hypothesis on $f$, the first contact is
at $s=0$ and the lemma is proved.
\end{proof}

These estimates have a lot of important consequences in our study of the
uniqueness. First, it gives us a technical lemma.

\begin{lem}\label{minmax}
Let $f:\R\rightarrow\R$ be a continuous function. We consider one of the two
cases below:
\begin{enumerate}
\item $\Ome=\R\times (-\frac{1}{2H},\frac{1}{2H})$ and $f$ is convex.
\item $\Ome=\R \times (-h_t(H),h_t(H))$ and $f$ satisfies a uniform
  $\rho_t(H)$-circle under condition. 
\end{enumerate}
Let $\boD$ denote the set of all solutions $u$ of \eqref{cmc} on $\Ome$
with $\phi_f$ as boundary value. Let $u_1$ and $u_2$ be in $\boD$ then
there exist $v^+$ and $v^-$ in $\boD$ such that: 
\begin{align*}
v^+&\ge \max(u_1,u_2)& v^-&\le \min(u_1,u_2)
\end{align*}
\end{lem}
\begin{proof}
Let $n$ be in $\N$; we define
$\Ome_n=\{(x,y)\in\Ome\,|\,-n-\sqrt{1/(2H)^2-y^2} \le x \le
n+\sqrt{1/(2H)^2-y^2}\}$. The boundary of $\Ome$ is composed of two
segments and two circle-arcs of curvature $2H$. 

Using Perron process (see \cite{CH,GT}), we can build
\begin{itemize}
\item the solution $v^+_n$ of \eqref{cmc} on $\Ome_n$ with
  $\max(u_1,u_2)$ on the boundary and
\item the solution $v^-_n$ of \eqref{cmc} on $\Ome_n$ with
  $\min(u_1,u_2)$ on the boundary.
\end{itemize}
To build $v^+_n$, we use subsolutions (let us observe that
$\max(u_1,u_2)$ is a subsolution). By maximum principle, every
subsolution is less than $w$ the solution of \eqref{esm} given by
Lemma \ref{supminimal}. Then, we can define $v^+_n$ as the supremum
over all subsolutions. $v^+_n$ takes the good boundary values on the
two segments because $\max(u_1,u_2)=w$ on it. For the two circle-arcs,
we use the barrier functions built by J.~Serrin in \cite{Se1}. For
$v^-_n$, we use supersolutions ($\min(u_1,u_2)$ is one). By maximum
principle, every supersolution satisfies to the under-bound of Lemmas
\ref{infconvex} or 
\ref{infcircle}. Then we define $v^-_n$ as the infimum of all
supersolutions. The half-circles and nodoids of Lemmas \ref{infconvex}
and \ref{infcircle} are used as barrier functions and give us the
boundary value of $v^-_n$ on the two segments. For the two circle
arcs, we use J.~Serrin arguments. 

On $\Ome_n$, we have $\max(u_1,u_2)\le v^+_n\le w$ then a subsequence
converges to $v^+$ on $\Ome$ and $v^+\in \boD$. Clearly $\max(u_1,u_2)\le
v^+$. The sequence $v^-_n$ is 
upper-bounded by $\min(u_1,u_2)$ and satisfies the under-bounds of
Lemmas \ref{infconvex} or \ref{infcircle}. Then a subsequence
converges to $v^-$ a solution of \eqref{cmc} with $\phi_f$ as boundary
value. Besides, $\min(u_1,u_2)\ge v^-$    
\end{proof}

With this Lemma, we can prove:

\begin{prop}\label{encadrement}
Let $f:\R\rightarrow\R$ be a continuous function. We consider one of the two
cases below:
\begin{enumerate}
\item $\Ome=\R\times (-\frac{1}{2H},\frac{1}{2H})$ and $f$ is convex.
\item $\Ome=\R \times (-h_t(H),h_t(H))$ and $f$ satisfies a uniform
  $\rho_t(H)$-circle under condition. 
\end{enumerate}
There exist $u_{max}$ and $u_{min}$ two solutions of \eqref{cmc} on $\Ome$
with $\phi_f$ as boundary value such that, for every solution $u$ of
\eqref{cmc} on $\Ome$ with $\phi_f$ as boundary value, we have:
$$u_{min}\le u\le u_{max}$$ 
\end{prop}
\begin{proof}
Let us denote $\boD$ the set of all solutions $u$ of \eqref{cmc} on
$\Ome$ with $\phi_f$ as boundary value; thanks to P.~Collin and R.~L\'opez,
$\boD$ is non-empty. So we define $u_{max}$ and $u_{min}$ at 
$p\in\Ome$ by: 
\begin{gather*}
u_{max}(p)=\sup_{u\in\boD}u(p)\\
u_{min}(p)=\inf_{u\in\boD}u(p)
\end{gather*}
By Lemma \ref{supminimal}, $u_{max}$ is well defined ; Lemmas
\ref{infconvex} and \ref{infcircle} ensure that $u_{min}>-\infty$. As
in the classiccal Perron process, it can be proved that $u_{max}$ and $u_{min}$
are two solutions of \eqref{cmc} on $\Ome$: in fact the argument we
need is that for every $u_1$ and $u_2$ in $\boD$ there exist
$u_3\in\boD$ that upper-bounds $\max(u_1,u_2)$ and $u_4\in\boD$ that
under-bounds $\min(u_1,u_2)$ (this is Lemma \ref{minmax}).

Using the solution $w$ of \eqref{esm} built in Lemma \ref{supminimal},
the half-cylinders of Lemma \ref{infconvex} or the nodoids of Lemma
\ref{infcircle} as barrier functions, we finally prove that $u_{max}$
and $u_{min}$ have $\phi_f$ as boundary value. Besides the construction gives
, for every $u\in\boD$:
$$u_{min}\le u\le u_{max}$$
\end{proof}

We have an important remark on these two solutions. For every
$(x,y)\in\Ome$, they satisfy :
\begin{gather}
u_{max}(x,y)=u_{max}(x,-y)\label{sym1}\\
u_{min}(x,y)=u_{min}(x,-y)\label{sym2}
\end{gather} 
This is due to the fact that both functions $(x,y)\mapsto
u_{max}(x,-y)$ and $(x,y)\mapsto u_{min}(x,-y)$ are in $\boD$.

\subsection{Upper-bounds}\label{remarque}

In this subsection, we look for explicit upper-bounds for solutions of
\eqref{cmc}. First, we have the following upper-bound:

\begin{prop}\label{upperbound}
Let $f:\R\rightarrow\R$ be a contiuous function and $x_0\in\R$. We assume that
$f$ is monotonous on $[x_0,+\infty)$. Let $u$ be a solution of \eqref{cmc} on
$\Ome=\R\times (-a,a)$ with $\phi_f$ as boundary value. Then for $x\ge
x_0+1/H$, we have:
$$
u(x,y)\le f(x)+\frac{1}{2H}
$$ 
\end{prop}
\begin{proof}
We only consider the case where $f$ is increasing on $[x_0,+\infty)$. We
consider $a\ge x_0+1/H$ and denote by $C((a-1/(2H),s),1/(2H))$ the horizontal
cylinder of axis $\{x=a-1/(2H)\}\cap \{z=s\}$ and radius $1/(2H)$. For big $s$
the cylinder $C((a-1/(2H),s),1/(2H))$ is above the graph of $u$. Let $s$
decrease until $s_0$ where the first contact happens. By maximum principle,
this first contact point is on the boundary at a point of first coordinate
$a'\in[a-1/(2H),a]$. We have $f(a')\ge s_0-1/(2H)$.

Since for every $s\ge s_0$, $C((a-1/(2H),s),1/(2H))$ is above the graph of $u$,
we have $u(a,y)\le s$. Then $u(a,y)\le s_0\le f(a')+1/(2H)$. Since $a'<a$ and
$f$ is increasing, $u(a,y)\le f(a) +1/(2H)$.
\end{proof}

Let us introduce a definition. Let $f:\R\rightarrow\R$ be a continous
  function. $f$ satifies \emph{a $R$-circle upper condition at $a\in\R$} if
  there exists in $\{(x,y)\in\R^2|\, y\ge f(x)\}$ a disk $D$ with radius $R$
  such that $(a,f(a))\in \partial D$.

\begin{rem*}
Let $D((a,s),R)$ denote the disk with center $(a,s)$ and radius $R$. For big
$s$, $D((a,s),R)$ is included in $\{(x,y)\in\R^2|\, y\ge f(x)\}$. Let $s$
decrases until the first contact with the graph of $f$, $f$ then satifies a
$R$-circle upper condition at the first coordinates of the contact points. In
changing $a$, we get all the abscissas where $f$ satifies a $R$-circle upper
condition. This implies that for every $a\in\R$ there is $a'\in[a-R,a+R]$
where $f$ satisfies a $R$-circle upper condition.
\end{rem*}
\begin{prop}\label{major}
Let $f:\R\rightarrow\R$ be a continuous function. Let $u$ be a solution of
\eqref{cmc} on $\Ome=\R\times (-l,l)$ with $\phi_f$ as boundary value. We
assume the $f$ satisfies a $1/(2H)$-circle upper condition at $x_0\in\R$. Then,
for every $y\in[-l,l]$, $u(x_0,y)\le f(x_0)$.
\end{prop}
\begin{proof}
Let $\Gamma(a,b)$ denote the circle of center $(a,b)$ and radius $1/(2H)$
which belongs to $\{z\ge f(x)\}$ and such that $(x_0,f(x_0))\in
\Gamma(a,b)$. Let us denote by $C((a,b+s),1/(2H))$ the horizontal cylinder of
axis $\{x=a\}\cap \{z=b+s\}$ and radius $1/(2H)$. For big $s$ 
the cylinder $C((a-1/(2H),s),1/(2H))$ is above the graph of $u$. Let $s$
decrease until the first contact happens. Because of maximum principle and the
existence of $\Gamma(a,b)$, this first contact happens for $s=0$. Then, on the
segment $I_{x_0}=\{x_0\}\times[-l,l]$, $u$ is upper-bounded by $f(x_0)$.
\end{proof}

Let $f:\R\rightarrow\R$ be a continuous function that satisfies a uniform
$R$-circle under condition. Let $a\in\R$  denote a point where $f$ satifies a
$R'$-circle upper condition. Since at $a$, there is a circle below and over
the graph of $f$, the graph of $f$ has a tangent. Then either $f'(a)$ exists or
$f'(a)=\pm\infty$. In all the cases, we can deal with the sign of the
derivative of $f$ at $a$. We then have a kind of Rolle's Theorem.

\begin{lem}\label{rolle}
Let $f:\R\rightarrow\R$ be a continuous function that satisfies a uniform
$R$-circle under condition. Let $a<b$ be two points where $f$ satifies a
$R'$-circle upper condition. We assume that $f'(a)>0$ and $f'(b)<0$. Then
there exists $c\in[a,b]$ such that:
\begin{enumerate}
\item $f$ satifies a $R'$-circle upper condition at $c$.
\item $f'(c)=0$.
\end{enumerate}
\end{lem}
\begin{proof}
Let $g$ denote the function defined by $g(x)=R'-\sqrt{R'^2-x^2}$ on
$[-R',R']$; its graph is a half-circle of radius $R'$. Since $f$ satifies a
$R'$-circle upper condition at $a$ and $f'(a)>0$, $f$ is upper bounded by
$f(a)+g(x-a)$ on $[a-R',a]$. In the same way, $f$ is upper-bounded 
by $f(b)+g(x-b)$ on $[b,b+R']$. Let $c\in[a,b]$ denote a point where
$f(c)=\max_{[a,b]}f$. Then $f(x)$ is upper-bounded by $m(x)$ on
$[a-R',b+R']$ where $m(x)$ is defined by:
$$
m(x)=
\begin{cases}
f(c)+g(x-a)& \text{ for }x\in[a-R',a]\\
f(c)&  \text{ for }x\in[a,b]\\
f(c)+g(x-b)& \text{ for }x\in[b,b+R']
\end{cases}
$$

This implies that $f$ satisfies a $R'$-circle upper condition at $c$ and then
$f'(c)=0$. 
\end{proof}


\section{The uniqueness of Collin's solutions}
\label{proofcollin}
The aim of this section is to prove the uniqueness of the solutions for
the Dirichlet problem studied by P. Collin in \cite{Co}. More precisely, we
have the following result.
\begin{thm}\label{uniciteconvex}
Let $f:\R\rightarrow\R$ be a convex function. Let $u$ and
$v$ be two solutions of \eqref{cmc} on $\Ome=\R\times
(-1/(2H),1/(2H))$ with $\phi_f$ as boundary value. Then $u=v$.
\end{thm}

The proof of Theorem \ref{uniciteconvex} is long, so the rest of the section
is devoted to it. In this proof, we shall use the differential $1$-form
$\omega_u$. If $u$ is a function on a domaine of $\R^2$, $\omega_u$ is defined
by :
$$
\omega_u=\frac{u_x}{\sqrt{1+|\nabla u|^2}}\dd y -
\frac{u_y}{\sqrt{1+|\nabla u|^2}}\dd x  
$$
with $u_x$ and $u_y$ the two first derivatives of $u$. When $u$ is a solution
of \eqref{cmc}, $\omega_u$ satisfies $\dd \omega_u=2H\dd x\wedge \dd y$ (see
\cite{Sp}). 
\subsection{Preliminaries}

By Proposition \ref{encadrement}, there are two solutions $u_{min}$ and
$u_{max}$ of \eqref{cmc} on $\Ome$ with $\phi_f$ as boundary value such that,
for every solution $u$ of the same Dirichlet problem, $u_{min}\le u\le
u_{max}$. Then to prove the uniqueness, it is sufficient to prove
: $u_{min}=u_{max}$. 

So let us assume that $u_{min}\neq u_{max}$; it is then known that
$u_{max}-u_{min}$ is unbounded on $\Ome$ \cite{CK}. By exchanging $x$ with
$-x$, we can assume that:
\begin{equation}\label{limit}
\lim_{x\rightarrow+\infty}\max_{I_x}(u_{max}-u_{min})=+\infty
\end{equation}
where $I_x=\{x\}\times[-1/(2H),1/(2H)]$. Let $c$ denote
$\max_{I_0}u_{max}-u_{min}$. Then there exists $\boD$ a connected component of
$\{u_{max}\ge u_{min}+2c\}$ that is included in
$\R_+\times[-1/(2H),1/(2H)]$. $\boD$ is unbounded. 

Since $f$ is convex, $f$ has a left derivative $f'_l$ and a right
derivative $f'_r$ at every point. These two functions increase and have the
same limit at 
$+\infty$. If $\lim_{+\infty}f'_l=\lim_{+\infty}f'_r<+\infty$, $f$ is
lipschitz continuous on $\R_+$. Then Theorem 5 in \cite{Ma2} contradicts
\eqref{limit}.

Then $f$ must satisfy
\begin{equation}\label{limitf'}
\lim_{+\infty}f'_l=\lim_{+\infty}f'_r=+\infty
\end{equation}

\subsection{Asymptotic behaviour of $u_{min}$}

Let $(x_n)$ be a real sequence with $\lim x_n=+\infty$. Let us define
$u_n$ on $\Ome$ by $u_n(x,y)=u_{min}(x+x_n,y)$. For $a\in\R$, let us denote by
$C^+(a)$ the circle arc $\{x\ge a\}\cap \{(x-a))^2+y^2=1/(4H^2)\}$. This
circle-arc has $(a,-1/(2H))$ and $(a,1/(2H))$ as end-points. Besides $C^+(a)$
contains the point $(a+1/(2H),0)$. We then have the following result. 
\begin{lem}\label{asymptconvex}
There exists $(x_n)$ a real increasing sequence with $\lim x_n=+\infty$ such
that $(u_n)$ has $C^+(0)$ as line of divergence.
\end{lem}

Before the proof, let us recall what is a line of divergence. We refer to
\cite{Ma1} for the details. Let $(v_n)$ be a sequence of solutions of
\eqref{cmc} and $N_n$ denote the upward pointing normal to the graph of
$v_n$. Let us assume that $N_n(P)$ tends to a horizontal unit vector $(\nu,0)$
($\nu\in\S^1$). Let $C$ denote the circle-arc in the $xy$-plane with radius
$1/(2H)$ such that $P\in C$ and $2H\nu$ is the curvature vector of $C$ at
$P$. $C$ is then a line of divergence of the sequence $(v_n)$. Let us extend
the definition of $\nu$ along $C$ by $2H\nu(Q)$ is the curvature vector of $C$
at $Q\in C$ ($\nu(Q)\in\S^1$). Then, for every $Q\in C$, $N_n(Q)\rightarrow
(\nu(Q),0)$. This implies that for every $C'$ a subarc of $C$:
$$
\lim_{n\rightarrow+\infty}\int_{C'}\omega_u=\ell(C')
$$ 
with $\ell(C')$ the length of $C'$. $C'$ is oriented such that $\nu$ is
left-hand side pointing along $C'$.
\begin{proof}[Proof of Lemma \ref{asymptconvex}]
Let $v_n$ be defined on $\Ome$ by $v_n(x,y)=u_{min}(x+n,y)$. The boundary
value of $v_n$ is $\phi_{f_n}$ with $f_n(x)=f(x+n)$. Because of
\eqref{limitf'}, $f_n$ is increasing on 
$[-1/H,+\infty)$ for big $n$. Then, by Proposition \ref{upperbound},
$v_n(0,0)\le f_n(0)+1/(2H)$. Now, let $\theta_n\in [0,\pi/2)$ such that
${f_n}'_l(0)\le \tan\theta_n\le {f_n}'_r(0)$. By Lemma \ref{infconvex}:
$$
v_n(\frac{1}{H},0)\ge -\frac{1}{\cos\theta_n} \sqrt{\frac{1}{4H^2}}+
\frac{1}{H}\tan\theta_n+ f_n(0)
$$
Because of \eqref{limitf'}, $\theta_n\rightarrow\pi/2$. Then:
\begin{equation*}
\begin{split}
v_n(\frac{1}{H},0)-v_n(0,0)&\ge \frac{1}{H\cos\theta_n}
(\sin\theta_n-\frac{1}{2})-\frac{1}{2H}\\
&\xrightarrow[n\rightarrow+\infty]{} +\infty 
\end{split}
\end{equation*}

Then the sequence of derivatives $\der{v_n}{x}$ can not stay upper-bounded on
$[0,1/H]\times\{0\}$. Then there exists a sequence $(a_n)$ in $[0,1/H]$ such
that: 
\begin{equation}\label{truc1}
\lim \der{v_n}{x}(a_n,0)=+\infty
\end{equation}
Let $x_n$ be defined by $n+a_n-1/(2H)$, we remak that $\lim x_n=+\infty$. We
consider $(u_n)$ the sequence of solution of \eqref{cmc} associated to
$(x_n)$. \eqref{truc1} becomes: 
$$
\lim \der{u_n}{x}(\frac{1}{2H},0)=+\infty
$$
Since $\der{u_n}{y}(1/(2H),0)=0$ by \eqref{sym2}, the limit normal to the
sequence of graphs over $(1/(2H),0)$ is $(-1,0,0)$. Then $C^+(0)$ is a line of
divergence for $(u_n)$. In considering a subsequence of $(x_n)$, we can assume
that it is increasing; this ends the proof.
\end{proof}

\subsection{End of Theorem \ref{uniciteconvex} proof}

Let $(x_n)$ be a sequence given by Lemma \ref{asymptconvex}. Let $\boD_n$
denote the following intersection: $$
\boD_n= \boD\bigcap\left\{(x,y)\in\Ome\,|\,
  x\le x_n+\sqrt{\frac{1}{4H^2}-y^2}\right\} 
$$
\begin{figure}
\begin{center}
\resizebox{0.8\linewidth}{!}{\input{figuniciteH3.pstex_t}}
\caption{\label{boN}}
\end{center}
\end{figure}
The boundary of $\boD_n$ is composed of $\partial\boD\cap \boD_n$ and
$\Gamma_n$ which is the part included in the circle-arc $C^+(x_n)$ (see Figure
\ref{boN}). Let $\widetilde{\omega}$ denote $\omega_{u_{max}}
-\omega_{u_{min}}$; we then have:  
$$
0=\int_{\partial\boD_n}\widetilde{\omega}=\int_{\partial\boD\cap \boD_n}
\widetilde{\omega}+\int_{\Gamma_n} \widetilde{\omega}
$$
Thanks to Lemma 2 in \cite{CK}, the integral on $\partial\boD\cap
\boD_n$ is negative; besides, since $(x_n)$ is increasing, it
decreases when $n$ is increasing. Besides we have:  
$$
0<-\int_{\partial\boD\cap \boD_n} \widetilde{\omega}=\int_{\Gamma_n}
\widetilde{\omega}\le 2\ell(\Gamma_n) 
$$
where $\ell(\Gamma_n)$ denote the length of $\Gamma_n$. Then $\ell(\Gamma_n)$
is far from $0$ uniformaly under-bounded. Because of Lemma \ref{asymptconvex}
and since $\Gamma_n\subset C^+(x_n)$, there exists $(\alpha_n)$ a sequence in
$[0,1]$ such that $\lim \alpha_n=1$ and
$$
\int_{\Gamma_n}\omega_{u_{min}}\ge \alpha_n\ell(\Gamma_n)
$$
Finally, for $n\ge n_0>0$, we have:
\begin{equation*}
\begin{split}
0<-\int_{\partial\boD\cap \boD_{n_0}}
\widetilde{\omega}\le-\int_{\partial\boD\cap \boD_n} \widetilde{\omega} 
&=\int_{\Gamma_n}\omega_{u_{max}}-\int_{\Gamma_n}\omega_{u_{min}}\\
&\le \ell(\Gamma_n)-\alpha_n\ell(\Gamma_n)\\
&\le (1-\alpha_n)\ell(\Gamma_n)\xrightarrow[n\rightarrow +\infty]{}0
\end{split}
\end{equation*}
Then we have a contradiction and Theorem \ref{uniciteconvex} is proved. \qed


\section{The uniqueness of L\'opez's solutions}
\label{prooflopez}
In this section, we prove the uniqueness of the solutions for
the Dirichlet problem studied by R. L\'opez in \cite{Lo1}. More precisely, we
have the following theorem.
\begin{thm}\label{unicitecircle}
Let $f:\R\rightarrow\R$ be a continuous function that satisfies a uniform
$\rho_t(H)$-circle under condition. Let $u$ and
$v$ be two solutions of \eqref{cmc} on $\Ome=\R\times
(-h_t(H),h_t(H))$ with $\phi_f$ as boundary value. Then $u=v$.
\end{thm}

The proof of this theorem is very similar to the one of Theorem
\ref{uniciteconvex}. The following of the section is devoted to it.




\subsection{Preliminaries}

By Proposition \ref{encadrement}, there are two solutions $u_{min}$ and
$u_{max}$ of \eqref{cmc} on $\Ome$ with $\phi_f$ as boundary value such that,
for every solution $u$ of the same Dirichlet problem, $u_{min}\le u\le
u_{max}$. Then to prove the uniqueness, it is sufficient to prove
: $u_{min}=u_{max}$. 

So let us assume that $u_{min}\neq u_{max}$; it is then known that
$u_{max}-u_{min}$ is unbounded on $\Ome$. By exchanging $x$ with $-x$, we can
assume that:
\begin{equation}\label{limitcircle}
\lim_{x\rightarrow+\infty}\max_{I_x}(u_{max}-u_{min})=+\infty
\end{equation}
where $I_x=\{x\}\times[-h_t(H),h_t(H)]$. Let $c$ denote
$\max_{I_0}u_{max}-u_{min}$. Then there exists $\boD$ a connected component of
$\{u_{max}\ge u_{min}+2c\}$ that is included in
$\R_+\times[-h_t(H),h_t(H)]$. $\boD$ is unbounded. 

Equation \eqref{limitcircle} has consequences. First, the existence of two
different solutions implies:
\begin{lem}
There exists $x_0\in\R_+$ such that $f$ is monotonous on $[x_0,+\infty)$.
\end{lem}
\begin{proof}
Let us consider the set $\boS$ of points where $f$ satisfies a $1/(2H)$-circle
upper condition. From a remark in Section \ref{remarque}, $\boS$ is non-empty
and is unbounded. Let us recall that, for every point in $\boS$, we can deal
with the sign of the derivative of $f$. First, we prove that there exists
$x_1\in\R^+$ such that, for every $x\in \boS\cap [x_1,+\infty)$, the sign of
$f'(x)$ is constant. If it is not true there exists two sequences $(a_n)$ and
$(b_n)$ in $\boS$ such that: 
\begin{itemize}
\item $\lim_{n\rightarrow +\infty}a_n=+\infty$ and $\lim_{n \rightarrow+\infty}
  b_n=+\infty$
\item $a_1<b_1<a_2<b_2<\cdots<a_n<b_n<\cdots$
\item for every $n$, $f'(a_n)$ is positive and $f'(b_n)$ is negative.
\end{itemize}
Because of Lemma \ref{rolle}, there exists a sequence $(c_n)$ in $\boS$
such that $a_n<c_n<b_n$ and $f'(c_n)=0$. Let $u$ denote a solution of
\eqref{cmc} on $\Ome$ with $\phi_f$ as boundary value. By Proposition
\ref{major}, $\max_{I_{c_n}}u\le f(c_n)$. Besides, since $f'(c_n)=0$,
Lemma \ref{infcircle} implies that $\min_{I_{c_n}}u\ge f(c_n)-(\rho_t(H)-t)$. 

So this implies that
$$
\max_{I_{c_n}}(u_{max}-u_{min})\le (\rho_t(H)-t)
$$
Since $\lim c_n=+\infty$, this contradicts \eqref{limitcircle}. Then there
exists $x_1\in\R^+$ such that, for every $x\in \boS\cap [x_1,+\infty)$, the
sign of $f'(x)$ is constant. We assume in the following that these derivatives
are positive. 

If there is no $x_0$ such that $f$ increases on $[x_0,+\infty)$, there is a
sequence $a_n\in [x_1,+\infty)$ such that:
\begin{itemize}
\item $(a_n)$ increases and $\lim a_n=+\infty$
\item for every $a_n$, $f(a_n)$ is a local maximum of $f$.
\end{itemize}
Since $f$ satisfies a $\rho_t(H)$-circle under condition, we remark that $f$
is differentiable at every $a_n$. Let $\Gamma((a_n-1/(2H),s)$ denote the
circle of center $(a_n-1/(2H),s)$ and radius $1/(2H)$. For big $s$,
$\Gamma((a_n-1/(2H),s)$ is above the graph of $f$. Let $s$ decrease until
$s_0$ where the first contact happens. We get a point $x$ where $f$ satisfies
a $1/(2H)$-circle upper condition. By what we proved above, $f'(x)>0$ then $x$
belongs to $[a_n-1/(2H),a_n]$. Let $b_n\in[x,a_n]$ denote a point where
$f(b_n)=\max_{[x,a_n]}f$. Since $f'(x)> 0$, $b_n\in(x,a_n]$ then
$f'(b_n)=0$. Using horizontal cylinders with $\Gamma((a_n-1/(2H),s)$ as
vertical section, we prove that $\max_{I_{b_n}}u\le f(x)+1/(2H)\le
f(b_n)+1/(2H)$ with $u$ a solution of \eqref{cmc} on $\Ome$ with $\phi_f$ as
boundary value. Besides since $f'(b_n)=0$, $\min_{I_{b_n}}u\ge
f(b_n)-(\rho_t(H)-t)$. This implies that:
$$
\max_{I_{b_n}}(u_{max}-u_{min})\le 1/(2H)+(\rho_t(H)-t)
$$ 
As $\lim b_n=+\infty$, the above inequation contradicts
\eqref{limitcircle}. The lemma is then proved.
\end{proof}

As in the above proof, we assume in the following of Theorem
\ref{unicitecircle} proof that $f$ is increasing on some $[x_0,+\infty)$. If
$f$ decreases the argument are similar to the one we are going to give.

From Theorem 5 in \cite{Ma2}, we know that $f(x+4/H)-f(x)$ can not stay
bounded when $x$ goes to $+\infty$. We even know that:
\begin{equation}\label{limitf-f}
\lim_{x\rightarrow+\infty} f(x+4/H)-f(x)=+\infty
\end{equation}

In this proof, this indentity plays the same role as \eqref{limitf'} in the
proof of Theorem \ref{uniciteconvex}.

\subsection{The asymptotic behaviour of $u_{min}$}

Let $(x_n)$ be a real sequence with $\lim x_n=+\infty$. Let us define
$u_n$ on $\Ome$ by $u_n(x,y)=u_{min}(x+x_n,y)$. For $a\in\R$, let us denote by
$C^+(a)$ the circle arc:
$$
\left\{\left(x-\left(a-\sqrt{\frac{1}{4H^2}-h_t(H)^2}\right)\right)^2+ y^2=
  \frac{1}{4H^2}\right\} \bigcap \left\{x\ge a\right\}
$$
This circle-arc has $(a,-h_t(H))$ and $(a,+h_t(H))$ as end-points. Besides
$C^+(a)$ contains the point $(a+K,0)$ with
$K=\dis\frac{1}{2H}-\sqrt{\frac{1}{4H^2}-h_t(H)^2}$. We then have the
following result. 
\begin{lem}\label{asymptcircle}
There exists $(x_n)$ a real increasing sequence with $\lim x_n=+\infty$ such
that $(u_n)$ has $C^+(0)$ as line of divergence.
\end{lem}

\begin{proof}
Let $v_n$ be defined on $\Ome$ by $v_n(x,y)=u_{min}(x+n,y)$. The boundary
value of $v_n$ is $\phi_{f_n}$ with $f_n(x)=f(x+n)$. For $n$ big enough, $f_n$
is increasing on 
$[1/H,+\infty)$; so, using Proposition \ref{upperbound}, $v_n(0,0)\le
f_n(0)+1/(2H)$. Now let us apply Lemma \ref{infcircle}, we get that
$v_n(4/H+\rho_t(H),0)\ge f_n(4/H)-(\rho_t(H)-t)$. To get this
under-bound, Lemma \ref{infcircle} is applied at $4/H$; the graph of $v_n$ is
then above a nodoid $\boN_t$ with horizontal axis in the vertical plane
$x=4/H+A$ ($0\le A\le \rho_t(H)$ since $f$ increases). Since $\boN_t$ is below
the graph $v_n(4/H+A,0)\ge f_n(4/H)-(\rho_t(H)-t)$ (see Figure
\ref{nodoid}). Now let us translate $\boN_t$ by the horizontal vector
$e_x=(1,0,0)$; since $f_n$ is increasing, the nodoid $\boN_t+se_x$ stays under
the graph since it does not cross its boundary. Then for $s=\rho_t(H)-A$ we
get $v_n(4/H+\rho_t(H),0)\ge f_n(4/H)-(\rho_t(H)-t)$. 
\begin{figure}
\begin{center}
\resizebox{0.39\linewidth}{!}{\input{figuniciteH4.pstex_t}}
\caption{\label{nodoid}}
\end{center}
\end{figure}

Then we have:
$$
v_n(4/H+\rho_t(H),0) -v_n(0,0)\ge f_n(4/H)-f_n(0)-\frac{1}{2H}-\rho_t(H)+t
$$

By \eqref{limitf-f}, $\lim v_n(4/H+\rho_t(H),0) -v_n(0,0)=+\infty$. Then the
sequence of derivatives $\der{v_n}{x}$ can not stay upper-bounded on
$[0,4/H+\rho_t(H)]\times\{0\}$. Then there exists a sequence $(a_n)$ in
$[0,4/H+\rho_t(H)]$ such that:
\begin{equation}\label{truc}
\lim \der{v_n}{x}(a_n,0)=+\infty
\end{equation}

Let us recall that $K$ denote
$\dis\frac{1}{2H}-\sqrt{\frac{1}{4H^2}-h_t(H)^2}$. Let $x_n$ be defined by
$n+a_n-K$, we remak that $\lim x_n=+\infty$. We consider $(u_n)$ the sequence of
solution of \eqref{cmc} associated to $(x_n)$. \eqref{truc} 
becomes:
$$
\lim \der{u_n}{x}(K,0)=+\infty
$$
Since $\der{u_n}{y}(K,0)=0$ by \eqref{sym2}, the limiting normal to the
sequence of graphs over $(K,0)$ is $(-1,0,0)$. Then $C^+(0)$ is a line of
divergence for $(u_n)$. In considering a subsequence of $(x_n)$, we can assume
that it is increasing; this ends the proof.
\end{proof}

\subsection{End of Theorem \ref{unicitecircle} proof}

Let $(x_n)$ be the sequence given by Lemma \ref{asymptcircle}. Let $\boD_n$
denote the following intersection: 
$$
\boD_n=\boD\cap\left\{(x,y)\in\Ome\,|\, x\le x_n+\sqrt{\frac{1}{4H^2}-y^2}-
  \sqrt{\frac{1}{4H^2}-h_t(H)^2}\right\} 
$$
The boundary of $\boD_n$ is composed of $\partial\boD\cap \boD_n$ and
$\Gamma_n$ which is the part included in the circle-arc $C^+(x_n)$. Let $
\widetilde{\omega}$ denote $\omega_{u_{max}} -\omega_{u_{min}}$; we then
have: 
$$
0=\int_{\partial\boD_n}\widetilde{\omega}=\int_{\partial\boD\cap \boD_n}
\widetilde{\omega}+\int_{\Gamma_n} \widetilde{\omega}
$$
On $\partial\boD\cap \boD_n$, the integral is negative; besides, since $(x_n)$
is increasing, it decreases when $n$ is increasing (Lemma 2 in
\cite{CK}). Besides, we have:
$$
0<-\int_{\partial\boD\cap \boD_n} \widetilde{\omega}=\int_{\Gamma_n}
\widetilde{\omega}\le 2\ell(\Gamma_n)
$$
where $\ell(\Gamma_n)$ denote the length of $\Gamma_n$. Then $\ell(\Gamma_n)$
is far from $0$ uniformaly under-bounded. Because of Lemma \ref{asymptcircle}
and since $\Gamma_n\subset C^+(x_n)$, there exists $(\alpha_n)$ a sequence in
$[0,1]$ such that $\lim \alpha_n=1$ and
$$
\int_{\Gamma_n} \omega_{u_{min}}\ge \alpha_n\ell(\Gamma_n)
$$
Finally, for $n\ge n_0>0$, we have:
\begin{equation*}
\begin{split}
0<-\int_{\partial\boD\cap \boD_{n_0}} \widetilde{\omega}\le
-\int_{\partial\boD\cap \boD_n} \widetilde{\omega} 
&=\int_{\Gamma_n} \omega_{u_{max}}-\int_{\Gamma_n} \omega_{u_{min}}\\
&\le \ell(\Gamma_n)-\alpha_n\ell(\Gamma_n)\\
&\le (1-\alpha_n)\ell(\Gamma_n)\xrightarrow[n\rightarrow +\infty]{}0
\end{split}
\end{equation*}
Then we have a contradiction and Theorem \ref{unicitecircle} is proved. \qed

Let us explain what are the differences if we assume that $f$ is decreasing
and not increasing. In this case, we have to study the asymptotic behaviour of
$u_{max}$. We prove that there exists a sequence $(x_n)$ with $\lim
x_n=+\infty$ such that $C^-(0)$ is line of divergence of $(u_n)$. Here $u_n$
is defined by $u_n(x,y)=u_{max}(x+x_n,y)$ and $C^-(a)$ denotes the circle-arc:
$$
\left\{\left(x-\left(a+\sqrt{\frac{1}{4H^2}-h_t(H)^2}\right)\right)^2+ y^2=
  \frac{1}{4H^2}\right\} \bigcap \left\{x\le a\right\}
$$
With this result, we can make the computations of the end of the proof.


\bigskip

\noindent Laurent Mazet

\noindent Universit\'e Paul Sabatier, MIG

\noindent Laboratoire Emile Picard. UMR 5580

\noindent 31062 Toulouse cedex 9, France.

\noindent E-mail: mazet@picard.ups-tlse.fr

\end{document}